\documentclass[12pt]{article}

\usepackage{amsfonts,pdfsync}
\usepackage{amsthm}

\usepackage{amsmath}
\usepackage{graphicx}
\usepackage{comment}
\usepackage[usenames,dvipsnames]{xcolor}
\usepackage{enumerate}
\usepackage{extarrows}
\usepackage{hyperref}

\usepackage{stackrel}
\usepackage{centernot}
\usepackage{caption}
\usepackage{subcaption}
\usepackage{fancyhdr}

\let\quoteOLD\quote
\def\quote{\quoteOLD\small}

\definecolor{labelkey}{cmyk}{0,0.8,1,0.5}
\definecolor{refkey}{cmyk}{0,0.8,1,0.5}

\newtheorem{theorem}{Theorem}[section]
\newtheorem{example0}{\sc Example}[subsection]

\newtheorem{proposition}{Proposition}

\newtheorem{lemma}{Lemma}
\newtheorem{remark}{Remark}

\numberwithin{equation}{section}
\numberwithin{theorem}{section}
\numberwithin{corollary}{section}
\numberwithin{proposition}{section}
\numberwithin{lemma}{section}
\numberwithin{definition}{section}
\numberwithin{remark}{section}

\makeatletter
\def\th@newremark{\th@remark\thm@headfont{\bfseries}}
\makeatletter

\def\boxit#1{\vbox{\hrule\hbox{\vrule\kern6pt
          \vbox{\kern6pt#1\kern6pt}\kern6pt\vrule}\hrule}}

\newcommand{\sid}[1]{{\color{black} #1}}

\newcommand{\Fbar}{\overline{F}}
\newcommand{\Hbar}{\overline{H}}

\newcommand{\Phibar}{\overline{\Phi}}

\newcommand{\Gbar}{\overline{G}}
\newcommand{\Lbar}{\overline{L}}

\newcommand{\R}{\Bbb{R}}

\newcommand{\rmd}{{\rm d}}

\newcommand{\halmos}{\quad\hfill\mbox{$\Box$}}

\newcommand{\wt}{\widetilde}

\newcommand{\beq}{\begin{equation}}
\newcommand{\eeq}{\end{equation}}

\newcommand{\be}{\begin{equation}}
\newcommand{\ee}{\end{equation}}
\newcommand{\bea}{\begin{eqnarray}}
\newcommand{\eea}{\end{eqnarray}}
\newcommand{\bean}{\begin{eqnarray*}}
\newcommand{\eean}{\end{eqnarray*}}
\newcommand{\ben}{\begin{equation*}}
\newcommand{\een}{\end{equation*}}
\newcommand{\ba}{\begin{aligned}}
\newcommand{\ea}{\end{aligned}}

\def\nexto{\kern -0.54em}

\newcommand{\EE}{\textbf{\rm E}}

\def\Kless{K^u}
\def\Kgreater{K^c}

\begin{document}

\bibliographystyle{plain}
\title{Extremes of Censored and Uncensored Lifetimes in Survival Data} 


\author{Ross Maller$^a$ and   Sidney Resnick$^b$ 
\thanks{  
Email:  Ross.Maller@anu.edu.au;  sir1@cornell.edu
}}


\maketitle

\begin{abstract} 

The i.i.d. censoring model for survival analysis assumes two independent sequences of i.i.d. positive random variables,
$(T_i^*)_{1\le i\le n}$ and $(U_i)_{1\le i\le n}$.
The data 
consists of observations on the random sequence 
$(T_i)=(\min(T_i^*,U_i))$ together with  accompanying  censor indicators.
Values of $T_i$ with $T_i^*\le U_i$ are said to be uncensored, 
those with $T_i^*> U_i$ are censored.
We assume that the distributions of the $T_i^*$ and $U_i$ 
are in the domain of attraction of the Gumbel distribution 
and obtain the asymptotic distributions, as sample size $n\to\infty$,  of the maximum values of
the censored and uncensored lifetimes in the data, and of statistics related to them.
These enable us to examine questions concerning the possible existence of cured individuals in the population.
\end{abstract}

\section{Introduction}\label{intro}
In this paper we consider the i.i.d. censoring model in survival analysis, motivated by the fact that, in observed survival data, it is sometimes the case that the lifetimes of some of the longest-lived individuals in the sample are censored at the limit of follow-up time.
This can be taken as indicative of the existence in the population 
of a proportion of  ``cured'' individuals, or individuals ``immune'' to the event of interest (death of a patient, or recurrence of a disease, etc.) 
Consequently it is of interest to analyse the maximum values of
the censored and uncensored lifetimes in the data, and compare their magnitudes.
In the present paper we assume a realistic class of distributions 
for the survival and censoring distributions -- namely, 
those in the domain of attraction of the Gumbel distribution --
and obtain the joint asymptotic distribution of these maxima, and of statistics derived from them, 
and examine questions related to the existence of cured individuals in the population.

\subsection{The Data Model}\label{s2}
We assume a general independent censoring model with right censoring.
We have two independent sequences of i.i.d. positive random variables 
$(T_i^*)_{1\le i\le n}$ and $(U_i)_{1\le i\le n}$ having   
cumulative distribution functions (cdfs) 
$\wt F$  and $G$  on $[0,\infty)$. The data in a sample of size $n$ consists of observations on the random sequence of 
(possibly censored) survival times $T_i=\min(T_i^*,U_i)$, together with  accompanying  censor indicators.
%
The censoring distribution $G$ is assumed proper (total mass 1), but  the distribution $\wt F$ of the  $T_i^*$  is in general improper, with mass at infinity corresponding to cured individuals (who, formally, live forever). We assume it to be  of the form
\begin{equation}    \label{FandF0}
  \wt F(t)= pF(t),\ t\ge 0,
\end{equation}
where $0<p\le 1$    
and $F$ is the proper distribution of the ``susceptible" individuals.
 Only susceptibles can experience the event of interest and have an uncensored failure time. 

\begin{figure}[ht]
\centering
\centering
\includegraphics[width=12cm]{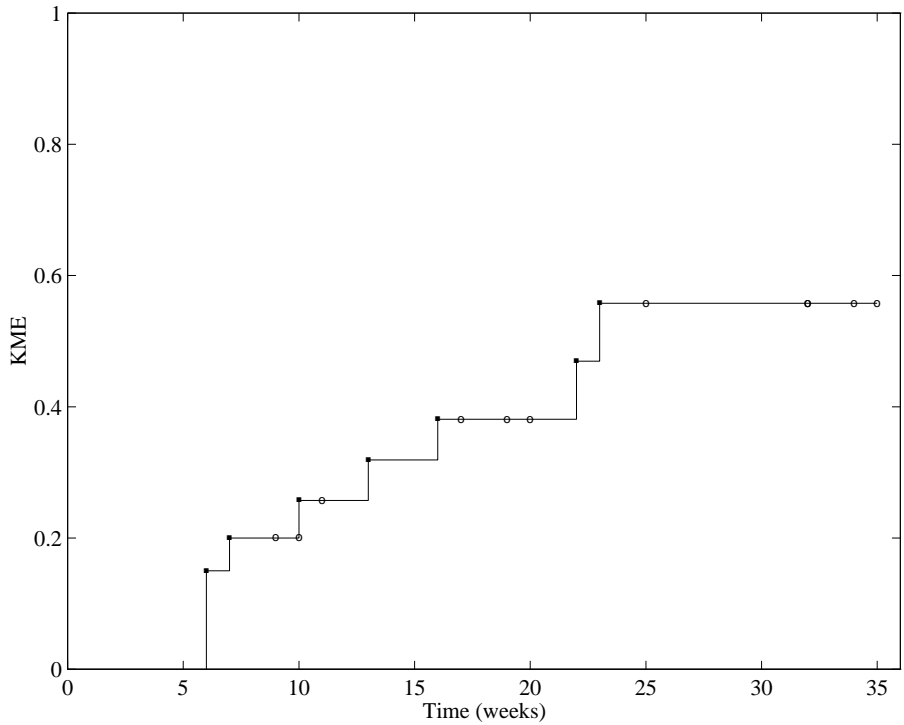}
\caption{Leukaemia Data}\label{fig:1}
\end{figure}

An informative way to display the sample  data is with the Kaplan-Meier estimator (KME) of the lifetime distribution; that is, the analogue of the empirical distribution function after censoring is taken into account. 
Figure 1 shows the KME constructed from data on  21 leukaemia patients (data from \cite{gehan:1965}, also in Figure 1.1, p.2,  of \cite{maller:zhoubook:1996}).\footnote{Figure \ref{fig:1} is a plot of a small, but real, data set. We include it as a schematic to display the features we are interested in. }
 The KME jumps at uncensored data times (full dots in Fig.1) 
and remains constant at censored points (open circles in Fig.1). 

%

A significant feature is the levelling of the KME below 1   at the right hand end (so the empirical distribution is improper) with a number of the largest observations being censored defining the level stretch. 
Such long-censored lifetimes  indicate the possibility of cured individuals being present in the population.
The useful information for this purpose is in the righthand end of the KME and of course these are the largest observations -- censored, in this case of interest -- suggesting an application of extreme value theory to study the  distribution of  the largest  {\it censored } lifetime.
Besides  this,  we  also want information on the 
largest {\it uncensored} lifetime, and, furthermore,  we need a comparison between the two.
The largest  uncensored lifetime in  Figure 1 is at 23 weeks, 
 the largest  censored  lifetime is at 35 weeks, and the 12 weeks difference between them is the length of the level stretch at the righthand end of the KME.

Our approach is to assume both distributions $F$ and $G$ 
are in the domain of attraction of the Gumbel distribution
and are comparable in terms of a certain balance condition on their hazard functions. 
 Such distributions  
 include the exponential, normal, lognormal, Weibull, and indeed most of the common distributions in use in survival analysis. 
 
 The following theorem   encapsulates our main findings.

\newpage

\begin{theorem}\label{thm:c>0}
Suppose $F$ and $G$ are both in the domain of attraction of the Gumbel distribution and their hazard functions satisfy a certain balance condition (Condition \eqref{e:assume} below)
depending on a parameter $\kappa\ge 0$.
Then for a sample of size $n$, we have the following results as $n\to\infty$.
\begin{enumerate}
  \item The largest uncensored lifetime and 
  the largest censored lifetime converge  jointly in distribution, after norming and centering,   
  to independent Gumbel random variables.
\item The largest uncensored lifetime and the largest lifetime (overall) converge   jointly, after norming and centering,    to a bivariate limiting random variable $(L_1,L_2)$.
\item The difference between the   largest \sid{observation}
  and  the largest \sid{un}censored lifetime converges  in distribution, after norming,   to the random variable $L:=L_2-L_1$, having cdf
 \be\label{Ldis}
 P[L\le x] =\frac{1}{1+\kappa e^{-x}}, \ x\ge 0.
 \ee
  \item The difference in Part 3, taken as a proportion of the largest observed lifetime, converges  in distribution (with no norming or centering needed),   to the random variable $R:=(L_2-L_1)/\max(L_1,L_2)$, having the distribution tail in \eqref{fin} below,
 depending only on the parameter $\kappa$. 
\end{enumerate} 
\end{theorem}  
The result in Part 3 of Theorem \ref{thm:c>0} is remarkably simple and explicit but its application in practice depends on knowing or estimating the norming sequence $a(n)$ in \eqref{e:defab}  below, as well as the parameter $\kappa$.
The result in Part 4 is more easily applicable, requiring only an estimate of  $\kappa$.
This parameter is related to the ratio of the hazard functions of the lifetime and  censoring distributions $F$ and $G$. We give further discussion of this, and examples, in an applications Section \ref{appl}. 

Another measure of the extent of followup in the sample is to count the number of censored lifetimes greater than the largest
  uncensored observation.
  In Section \ref{s3} we give the asymptotic distribution of this number,  under the same assumptions as in  Theorem \ref{thm:c>0}.

Theorem \ref{thm:c>0} will be proved  in Section \ref{pfth1}
and Theorem \ref{th:limCts}  in Section \ref{sec:count}. 
Our analysis prior to that, in Sections \ref{NOT} and \ref{pfp}, 
 proceeds by separating out, notionally,  the subsequences of censored and uncensored observations in the sample, applying extreme value techniques to each, then combining the results.

\section{Notation and Preliminary Results}\label{NOT}
Throughout we will assume both $F$ and $G$ are proper cdfs ($F(\infty)=G(\infty)=1$)
 with infinite right endpoints
 (the working can be modified to deal with finite right endpoints if they are the same for each distribution). 
Let $\Fbar(t)= 1-F(t)$ denote the survival function (tail function) of $F$, and similarly for  $\Gbar$. Let $H(t):=P(T_1\le t)$ 
be  the distribution of the observed survival times $(T_i)_{i\ge 1}$ with distribution tail $\Hbar(t)=1-H(t)=\Fbar(t)\Gbar(t)$. 


Relative to the sequence $\{(T_j^*,U_j), j \geq 1\}$, we define the random
indices $K_j^u$ and  $K_j^c$ by
\begin{align}\label{e:Kgreater}
  \Kless_0=&0, \quad \Kless_j =\inf \{m>\Kless_{j-1} : T_m^*\leq U_m\}, \ {\rm and} \cr
\Kgreater_0=&0, \quad \Kgreater_j =\inf \{m>\Kgreater_{j-1} : T_m^*>
                U_m\}.
\end{align}
Then the sequence $\{T_{\Kless_j}, j\ge 1\}$ selects out 
the subsequence of uncensored  observations in the sample,
and  the sequence $\{T_{\Kgreater_j}, j\ge 1\}$ selects out 
the subsequence of censored  observations.
Also  define
\begin{align} \label{e:ucount}
 & N_u(n) 
  =\sup\{m: \Kless_m \leq n\} \cr
  = &
  \{{\rm number\ of\ uncensored\ observations\ in\ a\ sample\ of\ size\ n}\}, 
\end{align}
and
\begin{align} \label{e:ccount}
  N_c(n) =n-N_u(n) 
  =
  \{{\rm number\ of\ censored\ observations\ in\ the\ sample\}}.
\end{align}

With the above notation the {\it largest uncensored lifetime}  in the first $n$ observations can be written as
  \be\label{mu}
  M_u(n): =  \max_{1\le i\le N_u(n)} T_{\Kless_i},
  \ee
and the {\it largest censored lifetime} is
$M_c(n):= \max_{1\le i\le N_c(n) }U_{\Kgreater_i}$.
The largest observation in the sample is  then
$$
M(n):= \max_{1\le i\le n}  T_i=
\max\big(M_u(n), M_c(n)\big).$$

\noindent \medskip {\it The D\'ecoupage de L\'evy.}\
A remarkable fact due to L\'evy (e.g., \cite[p.212]{resnickbook:2008})
is that, with $\Kless_j$ and $\Kgreater_j$ defined by \eqref{e:Kgreater}, 
 both  subsequences $\{(T_{\Kless_j},U_{\Kless_j} )\}$ and $\{(T_{\Kgreater_j},U_{\Kgreater_j}) \}$ are comprised of i.i.d. random vectors. 
Furthermore, the three sequences
\be\label{e:indep}
\{(T_{\Kless_j},U_{\Kless_j} ), j\geq
1\},\; \{(T_{\Kgreater_j},U_{\Kgreater_j}), j\geq 1\}, \;\{N_u(j),
j\geq 1\} 
\ee 
 are independent of each other, and the sequence $\{N_u(j)\}$ is a renewal counting function (a sum of i.i.d.  indicator variables rvs).
The distribution of the 2-vector 
$(T_{\Kless_1}, U_{\Kless_1})$ is the conditional distribution of
$(T_1,U_1)$  given $T_1^*\leq U_1$;
 that is,
$$
\big(T_{\Kless_1}, U_{\Kless_1}\big)
\stackrel{d}{=} 
\big((T_1,U_1)|T_1^*\leq U_1\big)
= \big((T_1^*,U_1)|T_1^*\leq U_1\big).
$$
We have for the distribution tail of an uncensored lifetime 
\begin{align}\label{e:tail<}
  P[  T_{\Kless_1}>x]=& P[T_1>x|T_1^*\leq U_1]
                         =\frac{P[U_1\geq T_1^*>x]}{ P[U_1\geq T_1^*] }
     =\frac{\int_x^\infty \bar G(s) F(\rmd s)}{\int_0^\infty \bar G(s)        F(\rmd s)}.
\end{align}
Likewise, for  a censored lifetime 
\begin{align}
\big(T_{\Kgreater_1}, U_{\Kless_1}\big)
\stackrel{d}{=} 
\big((T_1,U_1)|T_1^*> U_1\big)= \big((U_1,U_1)|T_1^*> U_1\big).
\end{align}
Interchanging $F$ and $G$ in 
\eqref{e:tail<}, we get the distribution tail of a censored lifetime as
\be\label{e:tail>}
  P[  T_{\Kgreater_1}>x]=
  P[T_1>x|T_1^*>U_1]=
  P[  U_{\Kgreater_1}>x]=\frac{\int_x^\infty \Fbar(s) G(\rmd s)}{\int_0^\infty \Fbar(s)  G(\rmd s)}.
  \ee

\noindent \medskip {\it The domain of attraction (DOA) of the Gumbel.}\
The Gumbel   distribution in standard form has  cdf 
  \be\label{e:gumbel}
  \Lambda (x)=\exp\{-e^{-x}\},\quad x\in \R.
  \ee
Throughout we will assume  both $F$ and $G$ are absolutely continuous and both are in the  domain of attraction of the Gumbel.
(Write this as $F\in D(\Lambda)$ and refer to 
\cite[Sect. 1.1]{resnickbook:2008} or 
\cite[Sect 1.2]{dehaan:ferreira:2006}
for background.)
This allows us to calculate asymptotic distributions of maxima of uncensored and censored lifetimes, using extreme value theory applicable to the Gumbel distribution.
In fact, consistent with the analysis in \cite{emvz:2020},  we will assume a little more than just the domain of attraction condition, namely,  that  both $F$ and $G$ are {\it  Von Mises distributions} \cite[p. 40]{resnickbook:2008} whose tail functions have  the form
  \be \label{e:barF}
    1-F(x)=\Fbar(x)=k_1\exp\Bigl \{-\int_{x_0}^x \frac{1}{f(u)}
             du\Bigr\},\quad x>x_0,
 \ee
 and
\be   \label{e:barG}
    1-G(x)=\Gbar(x)=k_2\exp \Bigl\{-\int_{x_0}^x \frac{1}{g(u)}
             du\Bigr\}, \quad x>x_0,      
\ee
  where  $f,g$ are absolutely   continuous  functions on $[x_0,\infty)$ 
 with densities $f',g'$ satisfying $f'(x)\to 0$, $g' (x)
  \to 0$, $x\to\infty$.  
  In \eqref{e:barF} and \eqref{e:barG} $x_0$ is a lower bound for the interval on which the representations hold and $k_1,k_2$ are positive constants. 

  An important result for us is that, under \eqref{e:barF} and \eqref{e:barG}, the product $\Fbar\times \Gbar$, which is the tail of     the distribution of the observed survival time $T_1=T_1^*\wedge U_1$,  
   is also  the tail of a Von
  Mises distribution, as  shown in the next theorem.
The proof of the theorem is in Section \ref{pfp}.
  \begin{theorem}[$\Hbar$ is the tail of a Von  Mises distribution]
  \label{prop:vonmises}
    If $\Fbar$ and $\Gbar$ are Von Mises
    distribution tails satisfying \eqref{e:barF} and \eqref{e:barG}, then
$\Hbar=\Fbar \times \Gbar$ is a Von Mises distribution tail with auxiliary
function $h:=fg/(f+g)$.
    \end{theorem}
  To  \eqref{e:barF} and \eqref{e:barG} we will add a third condition: 

\noindent {\bf Condition A.}\ 
  \be\label{e:assume}
  \lim_{x\to\infty}\frac{f(x)}{g(x)} =\kappa,\quad 0\leq \kappa<\infty.
  \ee
  \begin{remark}\label{rem1} 
{\rm 
  The  functions $f,g$ are called {\it auxiliary functions}; see \cite{resnickbook:2008}, p.26.
Differentiation of \eqref{e:barF} and \eqref{e:barG} shows that they are  the {\it reciprocal hazard functions} of $F$ and $G$ on the interval $[x_0,\infty)$. 
  Condition A specifies a certain kind of balance between the hazard functions corresponding to $F$ and $G$, and the  magnitude of $\kappa$  measures the relative heaviness of the  tails of $F$ and $G$.
      We discuss the practical implications of these facts in 
    Section \ref{appl}.  
      }
\end{remark}
 
The next step in our development is to compare the distribution tails of the censored and uncensored lifetimes with $\Hbar$, using the balance Condition A. 
Let
\begin{align} \label{e:pu}
  p_u =&   \{{\rm probability\ an\ observation\ is\ uncensored\}}
  =    P[T_i^*\leq U_i]
\end{align}
and 
\begin{align} \label{e:pc}
  p_c =&   \{{\rm probability\ an\ observation\ is\ censored\}}
  =   P[T_i^*> U_i]  =1-p_u.
\end{align}
A simple calculation gives the formulae
\ben 
p_u= \int_0^\infty \Gbar(s) F(\rmd s)
\quad {\rm and} \quad 
p_c= \int_0^\infty \Fbar(s) G(\rmd s).
\een

\begin{theorem}[Tail behaviour of the censored and uncensored lifetimes]
\label{prop:censortail}
  Suppose \eqref{e:barF}, \eqref{e:barG} and \eqref{e:assume}
  hold. Then the following are true.
  \begin{enumerate}
    \item There exists a non-decreasing function $U(x)$ such that
  \be\label{e:relateFG}
  \frac{1}{1- G(x)} =U \Bigl(   \frac{1}{1- F(x)} \Bigr). 
  \ee
 The  function $U(x)$ is regularly varying with index $\kappa\geq 0$ 
 (slowly varying when $\kappa=0$).
\item For all $\kappa\ge 0$, the tail of the uncensored lifetime distribution
(see   \eqref{e:tail<}) satisfies
  \be\label{e:uncenTail}
  P[  T_{\Kless_1}>x]
     \sim \frac{1}{(1+\kappa)p_u} \Hbar(x), \ x\to\infty.
     \ee
 When $\kappa>0$ the tail of the censored      lifetime distribution
     (see  \eqref{e:tail>})  satisfies
       \be\label{e:cenTail}
  P[ U_{\Kgreater_1}>x]
     \sim \frac{\kappa}{(1+\kappa)p_c}  \Hbar(x),\ x\to \infty.\ee
When $\kappa=0$  \eqref{e:cenTail} remains true in the sense that
 \be\label{e:c0}
 \lim_{x\to\infty}\frac{P[  U_{\Kgreater_1}>x]}{\Hbar(x)}
=\lim_{x\to\infty}\frac{P[  U_{\Kgreater_1}>x]}{P[T_{\Kless_1}>x]}
 =0.
 \ee
  \end{enumerate}
\end{theorem}
 \noindent In view of \eqref{e:tail<}, the relation \eqref{e:uncenTail} can be expressed as
  \be\label{e:uncenTail2}
\lim_{x\to\infty} \frac{1}{\Hbar(x)} 
  \int_x^\infty \Gbar(s) F(\rmd s)
= \frac{1}{1+\kappa}, 
     \ee
     valid for $0\le \kappa<\infty$. 
     The next theorem provides a partial converse to this.
\begin{theorem}[partial converse to \eqref{e:uncenTail2}]
\label{prop:prop3}
Suppose $\Fbar$ and $\Gbar$ are von Mises functions, thus 
satisfying \eqref{e:barF} and \eqref{e:barG}.
Assume that 
\be\label{12}
\lim_{x\to\infty} 
\frac{1}{\Hbar(x)} 
\int_{x}^{\infty} \Gbar (y)\rmd F(y)
 =k\in(0,1),
\ee
Then 
\be\label{0}
\lim_{x\to\infty} 
\frac{f(x)}{g(x)} =\frac{1-k}{k}.
\ee
\end{theorem}

 Theorem \ref{prop:prop3} does not cover the cases $k=0$ or $k=\infty$, but remains true in these cases under extra assumptions on $f$ or $g$. We omit details of this.

See Section \ref{pfp} for the proofs of  Theorems \ref{prop:vonmises}, \ref{prop:censortail} and \ref{prop:prop3}.
With these, we can complete the proof of Theorem \ref{thm:c>0}  in Section \ref{pfth1}.

\section[Number]{Numbers of Censored Lifetimes}\label{s3}
In this section, rather than the length of the level stretch at the right hand end of the KME, as in  Theorem \ref{thm:c>0}, we consider the number of censored observations that are bigger than the largest uncensored lifetime. A large number of such observations may be evidence for the presence of  immunes in the population. In this section we give the asymptotic distribution of this number. 

Recall the definitions of  
$\Kgreater _j$, 
$N_c(n)$ and $M_u(n)$ in \eqref{e:Kgreater}, \eqref{e:ccount} and \eqref{mu}.
We need also the number of censored lifetimes in the sample that exceed a value $t>0$, defined as
\be\label{e:cs}
N_c(>t):=\sum_{j=1}^{N_c(n)} 1_{[U_{\Kgreater _j}>t]}.
\ee

\begin{theorem}\label{th:limCts}
Assume  \eqref{e:barF}, \eqref{e:barG} and \eqref{e:assume}  and keep $0<\kappa <\infty$.

(i)\ Conditional on $M_u(n)$ and $ N_c(n)$, the
number $N_c(>M_u(n))$ is asymptotically, as $n\to\infty$, Poisson with parameter $\kappa E$ where $E$ is a unit exponential rv.  
By this we mean
\be\label{Ncon}
\lim_{n\to\infty} 
P\bigl[N_c(>M_u(n))=j\big|M_u(n), N_c(n)\bigr]
=P[ {\rm Poiss} (\kappa E) =j\, |  \, E], 
\ee
for $ \ j=0,1,2,\ldots$, 
where  Poiss$(\cdot)$ is a Poisson rv with the indicated parameter.

(ii)\ Unconditionally, $N_c(>M_u(n))$  is
asymptotically a geometric rv  with success probability
$p_\kappa:=\kappa/(1+\kappa)$ and mean $\kappa$;  
thus, 
\be\label{Ncon2}
\lim_{n\to\infty} 
P\bigl[N_c(>M_u(n))=j\bigr] = (1-p_\kappa) p_\kappa^{j}, \ j=0,1,\ldots.
\ee
\end{theorem}

\sid{The} proof of Theorem \ref{th:limCts} is in Section \ref{sec:count}. 
Before moving on to the proofs  we give some examples and applications.

\section{Examples and Applications}\label{appl}

 \subsection{Parameter $\kappa$ and the heaviness of the censoring}
The parameter $\kappa$ measures the relative heaviness of the censoring.  
When $\kappa=0$ in \eqref{e:assume},
 the hazard for $G$ is strongly dominated by the hazard for $F$, corresponding to relatively very light  censoring
(large values of $U$ are more likely than for $T^*$, so less censoring tends to occur).   Increasing values of $\kappa$ introduce progressively heavier censoring.

When $0<\kappa<\infty$ the hazards are comparable, asymptotically, but a finer classification is possible in terms of the tails of $F$ and $G$. 
Under \eqref{e:barF}, \eqref{e:barG} and \eqref{e:assume}
the function $U(x)$ in \eqref{e:relateFG} is regularly varying with index $\kappa\geq 0$, so we have
  (\cite[Proposition 0.8.(i), p.22]{resnickbook:2008}) 
  $$  \lim_{x\to\infty} \frac{U(x)}{x}
  =\begin{cases}
    \infty,& \text{ if }\kappa>1,\\
    0,& \text{ if }0\leq \kappa<1.
    \end{cases} 
    $$
  Therefore, from \eqref{e:relateFG},
\be\label{1}
\lim_{t\to\infty} \frac{\Fbar(t)}{\Gbar(t)}=\begin{cases}
    \infty,& \text{ if }\kappa>1,\\
    0,& \text{ if }0\leq \kappa<1.\end{cases} 
    \ee
So we see that the value $1$ for $\kappa$ is critical: for values
$0<\kappa<1$, the tail of $G$ dominates that of  $F$, \sid{censoring
variables tend to be bigger than lifetimes} and thus censoring
tends to be lighter; 
when  $\kappa>1$,  the tail of $F$ dominates that of  $G$ and
censoring tends to be heavier. 

 In a \sid{practical} situation when cured individuals are present in
 the population we expect  an intermediate value of $\kappa$  
and the observed value of the proportion $R$ in Part 4 of
Theorem \ref{thm:c>0}
 gives some information on this.
 A sample value of $ R$ close to $1$ means
\sid{the maximal uncensored lifetime is significantly smaller than the
  maximal censoring variable},
\sid{while  if an observed value of $ R$ is close  to $0$,
the maximal observation is approximately}
equal to \sid{the maximal uncensored
lifetime.}
\sid{Both situations are visible in a KME plot; in particular, when
  $R$ is close to $0$, the KME plot should be close to 1 at its right extreme.}
In this case there is little evidence of  cured individuals in the
population; i.e., \sid{$p \approx 1$} in \eqref{FandF0}. 
In a \sid{practical} situation we expect to see an intermediate value of $R$, with its distribution tail given by the expression in \eqref{fin} below. 

We see from \eqref{fin} that $P(R=0)=1$ if and only if $\kappa=0$, so 
a test of the hypothesis $H_0:\kappa=0$ serves as  a test for the
existence of a cure proportion.
Rejection of $H_0$ implies the \sid{possible}
 existence of a cure proportion, and evidence against $H_0$ is a
 sample value of $R$ significantly greater then 0. The
 distribution tail of $R$ in \eqref{fin} can be used to calculate
 critical values for the test if an estimate of $\kappa$ is available.  
 We follow up on these statistical issues elsewhere, and turn next to some examples of distributions to which the  theory applies. 

\subsection{Distributions in the DOA
of the Gumbel.}\label{DOAG}
\medskip\noindent{\it The Weibull distribution.}\ 
 The Weibull distribution is in the  domain of attraction of the Gumbel. We consider it in the form
 \be\label{web}
 F(x) = 1-e^{-\lambda x^{\alpha}},\ x\ge 0,
 \ee
 in which $\alpha>0$ is the shape parameter and $\lambda>0$ is the scale parameter. It has density
 $F'(x)= \lambda x^{\alpha-1}  e^{-\lambda x^{\alpha}}$, hazard function 
  $\lambda \alpha x^{\alpha-1}$, and the function $f(x)$ for \eqref{e:barF} is $f(x)=(\lambda\alpha)^{-1} x^{1-\alpha}$, taken for $x>x_0=1$, say.
    Since $\lim_{x\to\infty}f'(x)=0$, $F$ is in $D(\Lambda)$.
 Suppose $G$ is also Weibull with corresponding  parameters $\beta$ and $\mu$. Then  $g(x)$ for \eqref{e:barG} is $g(x)=(\beta\mu)^{-1} x^{1-\beta}$,
$x>1$, $G\in D(\Lambda)$,  and, as $x\to\infty$, 
 \ben
\frac{f(x)}{g(x)} =  \Big(\frac{\beta\mu}{\alpha\lambda}\Big) x^{\beta-\alpha} 
\to \kappa =
\begin{cases}
0, & \beta<\alpha;           \\
\mu/\lambda, &\beta=\alpha, \\
\infty, & \beta>\alpha.   
\end{cases} 
 \een
 Thus all 3 cases $\kappa=0$, $0<\kappa<\infty$, $\kappa=\infty$
 in   Theorems \ref{prop:censortail} and \ref{prop:prop3} can occur.

    \medskip\noindent{\it The exponential distribution.}\ 
This is the case $\alpha=\beta=1$ of the Weibull setup. 
Then  $f(x)=1/\lambda$, $g(x)=1/\mu$ and $\kappa=\mu/\lambda$.
   Thus  the censoring is  lighter or heavier according as $\mu<\lambda$        or vice-versa
   (the mean censoring variable is inversely proportional to $\mu$, so smaller values of $\mu$ give higher values of the censoring variable, hence lighter censoring).
Similar conclusions hold in the Weibull case.   
           
   \medskip\noindent{\it The lognormal distribution is in $D(\Lambda)$.}\   
We write the lognormal cdf in the form
\ben
F(x)= \frac{1}{\sqrt{2 \pi}} 
\int_{y=0}^x\exp\Big(-\frac{(\log y)^2}{2\sigma_F^2}\Big) \frac{\rmd y}{y}
=\Phi\Big(\frac{\log x}{\sigma_F}\Big),\ x>0,
   \een
where $\Phi(x)$ is the standard normal  cdf with tail $\Phibar(x)$ and density $\phi(x)$. 
Then $F\in D(\Lambda)$ (\cite[p.43]{resnickbook:2008}).
The reciprocal of the hazard function is
\ben
f(x)= 
\frac{x\sigma_F \Phibar\big(\log x/\sigma_F\big)}{\phi(\log x/\sigma_F)},\ x>0,
   \een
and since $\Phibar(z)\sim z^{-1} \phi(z) $ as $z\to\infty$, we have 
\ben
f(x) \sim
\frac{x\sigma_F^2}{\log x}, \ {\rm as}\ x\to\infty.
   \een
If the censoring distribution $G$ is lognormal with parameter $\sigma^2_G$, we have    
   \ben
   \lim_{x\to\infty} \frac{f(x)}{g(x)}= 
\frac{\sigma_F^2}{\sigma_G^2}, 
   \een
     Thus the censoring is heavier or lighter according as $\sigma_F>\sigma_G$ or vice-versa.  
      
   \medskip\noindent{\it The normal distribution.}\    
The normal  distribution is not usually used as a survival distribution, still we can consider a distribution with a normal-like tail and set 
$\Fbar(x)=\Phibar(x/\sigma_F)$,
 $\Gbar(x)=\Phibar(x/\sigma_G)$, $x\in \R$.
 Then  (\cite[p.43]{resnickbook:2008})
    $$
    f(x)\sim \frac{\sigma_F}{x},\quad g(x)\sim \frac{\sigma_G}{x},\quad
    x\to\infty,$$
    so
      \ben
  \lim_{x\to\infty}  \frac{f(x)}{g(x)} =
  \frac{\sigma_F}{\sigma_G}.
  \een
  Thus just as for the lognormal the censoring is heavier or lighter according as $\sigma_F>\sigma_G$ or vice-versa.

\medskip\noindent{\it The Weibull with lognormal censoring.}\ 
Suppose $F$ is Weibull with the cdf in \eqref{web}  
and $G$ is lognormal with 
\ben
g(x) \sim \frac{x\sigma_G^2}{\log x}, \ {\rm as}\ x\to\infty.
   \een
   Then 
   \ben
   \frac{f(x)}{g(x)} \sim
\frac{1}{\alpha\lambda \sigma_G^2}
\frac{\log x}{x^{\alpha}}
   \een
and the only possible case  for Theorems \ref{prop:censortail} is $\kappa=0$.

 Similar examples to the above can be constructed from the    gamma distribution which is in $D(\Lambda)$ (\cite[p.34]{dehaan:ferreira:2006}).

\subsection{Comments and related literature.}\label{lit}


A variety of models have been developed to analyse lifetime data
of the kind displayed in Figure \ref{fig:1}, 
in which  a proportion of the population 
may be long-term survivors.
A systematic formulation and treatment of these issues is in  
\cite{maller:zhoubook:1996} 
to which we refer for further background information.
Since 1996 there has been a steady increase in interest in cure models.
For more recent reviews, we mention 
\cite{oblc:2012}, 
\cite{pt:2014},  
\cite{ti:2014}  
 and \cite{av:2018}. 
An earlier paper along the lines of the present analysis is \cite{maller:zhou:1993}.
A recent paper also assuming a lifetime distribution in the domain of attraction of the Gumbel is 
\cite{emvz:2020}. 
 In \cite{ev:2018}  
a lifetime distribution in the domain of attraction of the Fr\'echet distribution is  assumed.

\section{Proofs of Theorems}\label{pfp}
    
    \medskip \noindent {\bf Proof of Theorem \ref{prop:vonmises}.} 
Assume \eqref{e:barF} and \eqref{e:barG} with the positive constants $k_1$, $k_2$ and $x_0$. Then with $k=k_1k_2$ and  $h=fg/(f+g)$, we have for $x\ge x_0$
  \begin{align*}
    \Fbar(x)\Gbar(x)=&k \exp\Bigl \{-\int_{x_0}^x
 \Bigl(\frac{1}{f(u)} + \frac{1}{g(u)} \Bigr) du\Bigr\}
  =k \exp\Bigl\{-\int_{x_0}^x\Bigl(\frac{1}{h(u)} \Bigr) \rmd u\Bigr\}.
  \end{align*}
Taking derivatives we get
  \begin{align*}
  h'=
\Bigl(\frac{fg}{f+g} \Bigr)' =\frac{f'g +g'f}{f+g}-
    \frac{fg}{(f+g)^2} (f'+g')=:A+B.
    \end{align*}
    Now
    $|A(x)  |\leq |f'(x)  |+|g' (x)| \to 0      $ and also
\bean
|B (x)|
&\leq &
 \frac{f(x)  g (x) }{f^2(x)  +g^2 (x) +2f (x)g(x) } 
(|f'(x)  |+|g' (x)|)\cr
& \leq&
 \frac 12  (|f' (x) |+|g' (x) |)\to 0,\ {\rm as}\ x\to\infty.
\eean
It follows that  $h' (x)   \to 0$, $x\to\infty$, 
and $h$ has the required property for $\Hbar$ to be  a Von Mises distribution tail.         \halmos

    \medskip \noindent {\bf Proof of Theorem \ref{prop:censortail}.} 
When $\kappa>0$ the assertion in \eqref{e:relateFG} follows from \cite[Theorem 2.1,
page 249]{dehaan:1974equiv}, and we only comment on the case $\kappa=0$
where $U$ must be slowly varying. Define the non-decreasing functions
$$U_F=\frac{1}{1-F} \quad {\rm and}\quad
U_G=\frac{1}{1-G},
$$
with inverse functions $U_F^\leftarrow $ and $U_G^\leftarrow$.
Then \eqref{e:barF} and \eqref{e:barG} imply
$$\lim_{t\to\infty} \frac{U_F(t+xf(t))}{U_F(t)}=e^x.$$
Inverting this we obtain 
$$\lim_{t\to\infty}\frac{ U_F^\leftarrow (tx)-U_F^\leftarrow
  (t)}{f(U_F^\leftarrow (t))} =\log x,\ x>0.$$
Comparable expressions hold for $U_G$ and $U_G^\leftarrow$. 
The assumption  $\kappa=0$ implies that there exists $\varepsilon (t) \to 0$ with $f(t)=\varepsilon (t) g(t).$ Define
$$U(x)=U_G\circ U_F^\leftarrow (x).$$
Then for $x>0$ (but fixed),
\begin{align*}
\lim_{t\to\infty} \frac{U(tx)}{U(t)}
=&
\lim_{t\to\infty} 
\frac{1}{U_G(U_F^\leftarrow (t))}\times 
  U_G\Bigl(\frac{U_F^\leftarrow (tx)-U_F^\leftarrow (t)}{f(U_F^\leftarrow
  (t)) }\cdot f(U_F^\leftarrow (t)) +U_F^\leftarrow (t)\Bigr)\cr
  &                \\
  =&\lim_{s\to\infty } \frac{U_G(\log x \cdot f(s) +s) }{U_G(s)}
     =\lim_{s\to\infty } \frac{U_G(\log x \cdot \varepsilon(s)g(s) +s) }{U_G(s)}.
\end{align*}
Given $\varepsilon>0$, for all large $s$, $|\varepsilon(s) \log x| \leq
\varepsilon$ so
$$\limsup_{s\to\infty}\frac{U(tx)}{U(t)} \leq e^{\varepsilon},$$
and similarly the $\liminf $ is bounded below by $e^{-\varepsilon}$. This
shows that $U$ is slowly varying when $\kappa=0$.

To prove  \eqref{e:uncenTail}, we have from \eqref{e:relateFG}
\ben
\int_x^\infty \Gbar(s)F(\rmd s)
=\int_x^\infty \frac{1}{U(1/\Fbar(v))} F(dv)
  =\int_{1/\Fbar(x)}^\infty \frac{1}{U(v)v^2} dv,
  \een
 and applying Karamata's theorem   (\cite[p. 17]{resnickbook:2008}),
 this is   asymptotic to
 \ben 
\Bigl(\frac{1}{1+\kappa} \Bigr) \Bigl(
         \frac{1-F(x)}{U(1/\Fbar(x))}\Bigr)
          = \Bigl(\frac{1}{1+\kappa} \Bigr)  \Fbar(x)\Gbar(x)
         = \Bigl(\frac{1}{1+\kappa} \Bigr) \Hbar(x), \ {\rm as}\ x\to\infty.
\een

When $\kappa>0$, \eqref{e:cenTail} is proved by interchanging the roles of
$F,G$ and $f,g$ 
(so $\kappa$ is replaced by  $1/\kappa$)
 and then applying the proof of
\eqref{e:uncenTail}. When $\kappa=0$, we still have $\int_x^\infty \Gbar(s)F(\rmd s) \sim \Fbar(x)\Gbar(x),\,x\to \infty$, and from Fubini's
theorem
$$\int_x^\infty \Fbar(s)G(\rmd s)=\Fbar(x)\Gbar(x)-\int_x^\infty 
\Gbar(s)F(\rmd s)$$
so
$$
\frac{\int_x^\infty \Fbar(s)G(\rmd s)}{\Fbar(x)\Gbar(x)}=1-\frac{\int_x^\infty \Gbar(s) F(\rmd s)}{\Fbar(x)\Gbar(x)} \to 1-1=0.
$$
This completes the proof of  Theorem \ref{prop:censortail}. \halmos

\bigskip\noindent {\bf Proof of  Theorem  \ref{prop:prop3}:}\
This will be a consequence of the following lemma. \sid{We need the
following concept: A function $f(x)$ is {\it self-neglecting} if it
is positive on $[x_0,\infty)$ for some $x_0$ and 
satisfies 
\be\label{4a}
\lim_{x\to\infty} 
\frac{f(x)}{f(  x+ y f(x)     )}  =1, \ y>0.
\ee
The convergence in \eqref{4a} is locally uniform in $y$
(\cite[p.41]{resnickbook:2008}, 
\cite[p.120, Sect. 2.11]{bingham:goldie:teugels:1987}).}

\begin{lemma}\label{lem1}
Suppose  $f(x)$ is  self-neglecting and  $H$ is any  distribution function satisfying
\be\label{2z}
\int_{x}^{\infty} \big(\Hbar(y)/f(y)\big) \rmd y   \sim k\Hbar(x), 
\ x\to\infty, 
\ee
where $0<k<1$.
Then $H$ is in the  domain of attraction of the Gumbel with an auxiliary function $h(x)$ satisfying $h(x) \sim kf(x)$ as $x\to\infty$. 
\end{lemma}

\noindent {\bf Proof of  Lemma  \ref{lem1}:}\
Assume \eqref{2z} with $f$ satisfying \eqref{4a} and define
\be\label{chidef}
\chi(x)= \frac{1}{ \Hbar(y)/f(y)}\int_{x}^{\infty} \big(\Hbar(y)/f(y)\big) \rmd y, \ x\ge x_0.
\ee
Then
\ben
\exp\Big(-\int_{x_0}^x \frac{1}{\chi(y)} \rmd y\Big)
= 
\frac{\int_{x}^{\infty} \big(\Hbar(y)/f(y)\big) \rmd y}
{\int_{x_0}^{\infty} \big(\Hbar(y)/f(y)\big) \rmd y}
=: \Lbar(x)
\een
is the tail of a cdf $L$ defined on $[x_0,\infty)$, 
and by \eqref{2z} 
\be\label{5.2}
\Lbar(x) \sim \frac{k \Hbar(x)}
{\int_{1}^{\infty} \big(\Hbar(y)/f(y)\big) \rmd y},\ {\rm as}\ x\to\infty.
\ee
From \eqref{2z} and \eqref{chidef} we have
\ben
\chi(x)= \frac{1}{ \Hbar(x)/f(x)}\int_{x}^{\infty} \big(\Hbar(y)/f(y)\big) \rmd y 
\sim
 \frac{k\Hbar(x)}{ \Hbar(x)/f(x)}
 =kf(x), 
 \een 
 \sid{and it follows that  $\chi$ is self-neglecting from the local
   uniform convergence in \eqref{4a}.}
Consequently $L$ is in the  domain of attraction of the Gumbel
with auxiliary function $\chi$, 
and since by \eqref{5.2} $\Lbar(x)$ is asymptotically equivalent to a constant times $\Hbar(x)$, also $H$  is in the  domain of attraction of the Gumbel with auxiliary function $h$, say.
Again since $\Lbar(x)$ is asymptotically equivalent to a constant
times $\Hbar(x)$, the auxiliary functions of $L$ and $H$ can be taken
the same  (\cite[p. 67, Proposition 1.19]{resnickbook:2008}).
Thus $h(x)\sim \chi(x)\sim kf(x)$ as $x\to\infty$.

 \halmos

To complete the  proof of  Theorem  \ref{prop:prop3}, assume \eqref{12} and set $\Hbar(x)=\Fbar(x)\Gbar(x)$ 
where  $\Fbar$ and $\Gbar$ satisfy \eqref{e:barF} and \eqref{e:barG}
with auxiliary functions $f$ and $g$.
Apply Lemma \ref{lem1} to get
\ben
h(x) = \frac{f(x)g(x)}{f(x)+g(x)} \sim kf(x).
\een
This implies
$f(x)\sim (k^{-1}-1) g(x)$ as required in \eqref{0}.  
\halmos

%
%

\section{Proof of Theorem \ref{thm:c>0}.}\label{pfth1}

Assume \eqref{e:barF}, \eqref{e:barG} and
\eqref{e:assume}.
Recall from Theorem \ref{prop:vonmises}
that $\Hbar$ is then the tail of a Von  Mises distribution with auxiliary function $h$. 
Thus  $H$ is in the domain of attraction of the  Gumbel   distribution with the cdf  $\Lambda(x)$ in \eqref{e:gumbel}.
Analytically, this means that its  auxiliary function $h$ has the property
\ben 
\lim_{t\to\infty}\frac{\Hbar(t+xh(t))}{\Hbar(t)}=e^{-x},\ x>0, 
\een
and any positive  sequences $a(n)$ and $b(n)$ satisfying
\be\label{e:defab} 
\lim_{n\to\infty} n\Hbar(b(n))
=1
\quad {\rm and}\quad
a(n)= h(b(n))    
\ee
also satisfy
\be\label{imp}
\lim_{n\to\infty}n\Hbar(a(n)x+b(n))=e^{-x}, \ x\in \R.
\ee
Hence  $a(n)$ and $b(n)$ are the appropriate norming and centering sequences such that
 the maximum observation has the Gumbel limit:
\be\label{e:gumbelLim}
\lim_{n\to\infty} P\Bigl[ \frac{M(n) -b(n)}{a(n)}
\leq x\Bigr] =\Lambda(x),\ x\in \R,
\ee
 see  \cite{resnickbook:2008}, Prop. 1.1, p.40.

Further to this, there is also functional convergence to an extremal process\footnote{For background on extremal processes, see  \cite{resnickbook:2008}, Sect. 4.3, p.179.
}
 with standard Gumbel marginals, $(Y(t))_{t>0}$.
By this we mean that \eqref{e:gumbelLim} can be extended to 
\be\label{e:fL}
 \Bigl(\frac{M(\lfloor nt\rfloor) -b(n  )}{a(n)} \Bigr)_{t>0}
\Rightarrow 
\bigl(Y(  t) \bigr)_{t>0},
\ee
where the convergence is in  $D(0,\infty)$, the space of right continuous
$\R$-valued functions with finite left limits on $(0,\infty)$, and $(Y(t))$ satisfies
\be\label{e:fLim}
P[Y(t)\leq x] =\Lambda^t(x)
=\exp\{-te^{-x}\},\ x\in \R.
\ee
See  \cite{resnickbook:2008}, Prop. 4.20, p.211.

We can prove an analogous result for the joint convergence of the  extremes $M_u(n)$ and  $M_c(n)$, after recalling  the independence property in  \eqref{e:indep} which means we may analyse them separately.
The limits will involve two independent  extremal processes  $\{Y_u(t),t>0 \}$   and $\{Y_c(t),t>0 \}$, each   with  standard  Gumbel   marginals, as in \eqref{e:fLim}.

First consider the limit distribution for uncensored lifetimes.
Recall from \eqref{mu} that
 $ M_u(n)$ is the maximum of $N_u(n)$ i.i.d. copies of 
 $T_{\Kless_1}$. 
To analyse this we look first at the maximum of $n$ i.i.d. copies of 
 $T_{\Kless_1}$. 
 A standard asymptotic using \eqref{imp} gives
\begin{align}\label{r1}
\lim_{n\to\infty}  P\Bigl[ \frac{\bigvee_{i=1}^n T_{\Kless_i} -b(n)}{a(n)}   \leq  x\Bigr]
  =&
  \lim_{n\to\infty}  \Bigl( 1-\frac{nP[T_{\Kless_1}>a(n)x+b(n)]}{n} \Bigr)^n\cr
=&\exp\{-\lim_{n\to\infty} nP[T_{\Kless_1}>a(n)x+b(n)]\}.
\end{align}
Use \eqref{e:relateFG} to write the RHS of \eqref{r1} as 
\begin{align*}
&\exp\Bigl\{ - \frac{1}{(1+\kappa) p_u }  \lim_{n\to\infty}  n\bar
   H(a(n)x+b(n))\Bigr\}
\end{align*}
and because of \eqref{e:gumbelLim} this equals
\be\label{r2}
 \big(\Lambda(x)\big) ^{((1+\kappa)p_u)^{-1}} 
  =  P\Bigl[ Y_u\Bigl(\frac{1}{p_u (1+\kappa) }\Bigr) \leq x \Bigr],\ x>0.
\ee

\subsection{The case $\kappa>0$.}\label{subsec:c>0}
Keeping $\kappa>0$ now, we have for censored lifetimes, similar to \eqref{r1} and \eqref{r2}, 
$$
\lim_{n\to\infty}
P\Bigl[  \frac{ \bigvee_{i=1}^n T_{\Kgreater_i}   -b(n)}{a(n)}    \leq x\Bigr]=
  P\Bigl[Y_c\Bigl(\frac{\kappa}{p_c (1+\kappa) } \Bigr) \leq
  x\Bigr], \ x>0.
  $$
We can combine the separate convergences  into a bivariate functional limit theorem  using \cite[Proposition 4.20, p.211]{resnickbook:2008} again, after noting that, 
by the weak law of large numbers, $ N_u(n)/n \to p_u $ and $N_c(n)/n \to p_c$.
 So we get for $t>0$,   $\kappa>0$, 
  \begin{align}
 \Bigl(\frac{  \bigvee_{i=1}^{\lfloor nt\rfloor}
  T_{\Kless_i} -b(n )}{a(n)}, \,
  & 
\frac{ \bigvee_{i=1}^{ \lfloor nt\rfloor} T_{\Kgreater_i} -b(n  )}{a(n)},\,
\frac{N_u(n)}{n}, \,
 \frac{N_c(n)}{n}\Bigr)
\nonumber\\
&
\Rightarrow 
\Bigl(Y_u\Bigl(\frac{1}{p_u(1+\kappa)  }  t\Bigr),\,
Y_c\Bigl(\frac{\kappa}{p_c (1+\kappa) } t\Bigr),\,
 p_u, \,  p_c\Bigr),\label{e:fclt}
  \end{align}
where the convergence is in $D((0,\infty)\mapsto \R^2)\times  [0,1]^2
\mapsto \R^2 \times  [0,1]^2$. 

Now apply the almost surely continuous scaling map $(x(\cdot),y(\cdot), a,b)\mapsto
(x(a),y(b))$ from 
$D((0,\infty)\mapsto \R^2)\times  [0,1]^2
\mapsto \R^2 $
 to \eqref{e:fclt} to deduce
\be \label{e:megilla}
\Bigl(  \frac{M_u(n) -b(n)}{a(n)},\,
  \frac{M_c(n) -b(n)}{a(n)}  \Bigr)
                                                         \Rightarrow
   \Bigl(Y_u\Bigl( \frac{1}{1+\kappa} \Bigr) ,
Y_c\Bigl(\frac{\kappa}{1+\kappa} \Bigr) \Bigr). 
\ee 
The first component on the left in \eqref{e:megilla} is the 
centered and normed 
maximal uncensored observation
and the second component is the centered and normed  maximal censored lifetime. 
It follows from \eqref{e:megilla} that
\begin{align}
&
     \Bigl(  \frac{M_u(n) -b(n)}{a(n)},  \frac{M(n) -b(n)}{a(n)}
      \Bigr)\nonumber \\
=&
\Bigl(   \frac{M_u(n) -b(n)}{a(n)},\
 \frac{  M_u(n)\vee  M_c(n)  -b(n)}{a(n)}
       \Bigr)\nonumber \\
\Rightarrow &
   \Bigl(Y_u\Bigl( \frac{1}{1+\kappa} \Bigr), \ 
Y_u\Bigl( \frac{1}{1+\kappa} \Bigr)\vee Y_c\Bigl(\frac{\kappa}{1+\kappa}
 \Bigr)\Bigr)                                                              \nonumber \\
  =:&
  (L_1,L_2).\label{e:L1L2}
\end{align}
The first component on the left in \eqref{e:L1L2}  is the 
centered and normed 
maximal uncensored lifetime,
while the second is the 
centered and normed 
maximal observation. On the right, the components are in terms of rescaled standard Gumbel extremal processes.
  Note that the components $(L_1,L_2) $ are not independent.
Using \eqref{e:fLim},
they satisfy
\begin{align}
P[L_1=L_2]=&
P\Big[Y_u \Bigl( \frac{1}{1+\kappa} \Bigr)  >
Y_c \Bigl( \frac{\kappa}{1+\kappa}\Bigl)  \Bigr]
\nonumber \\
=&\int_\R\Lambda^{\kappa/(1+\kappa)}(x)\Lambda^{1/(1+\kappa} (\rmd x)
  =\frac{1}{1+\kappa}.\label{e:LimEq}
\end{align}

Next we need the asymptotic distribution of the difference between the largest observed lifetime and the
largest uncensored lifetime. For this,  take differences in \eqref{e:L1L2}  to get 
\begin{align*}   
 \frac{ M_u(n) 
\bigvee M_c(n)-b(n)}{a(n)} 
 & -
 \frac{ M_u(n) -b(n)}{a(n)}
 \nonumber\\
&    \Rightarrow 
Y_u\Bigl( \frac{1}{1+\kappa} \Bigr) \vee Y_c\Bigl(\frac{\kappa}{1+\kappa} \Bigr) 
  -  Y_u\Bigl( \frac{1}{1+\kappa} \Bigr)   \cr
  & =
  L_2-L_1=: L.
  \end{align*}
  Note that it's important here that the centering sequence $b(n)$ is the same for both components.   
  Since for $x>0$
\be\label{yb}
   P \Bigl[Y_c\Bigl( \frac{\kappa}{1+\kappa} \Bigr) \vee y>x  \Bigr]
=
\begin{cases}
 P \Bigl[Y_c\Bigl( \frac{\kappa}{1+\kappa}\Bigr) >x\Bigr],& \text{ if } y<x,\\
    1,& \text{ if }y>x,
    \end{cases} 
\ee
we have, using \eqref{yb}, 
  \begin{align}
    P[L>x]=&
        P \Bigl[Y_u\Bigl( \frac{1}{1+\kappa} \Bigr) \vee Y_c\Bigl(\frac{\kappa}{1+\kappa}\Bigr)
  -  Y_u\Bigl( \frac{1}{1+\kappa}\Bigr) >x\Bigr) \Bigr]  \cr
         =&
\int_\R (1-\Lambda^{\kappa/(1+\kappa)}(x+y)\Lambda ^{1/(1+\kappa)} (\rmd y).
  \end{align}
  Evaluating this we obtain
  \be\label{e:Ltail}
     P[L>x]= \frac{\kappa}{e^x+\kappa}, \ x\ge 0, 
\ee
We add to this some mass at 0:  from \eqref{e:LimEq} we have
  $$P[L=0]=\frac{1}{1+\kappa}.$$
Thus we arrive at \eqref{Ldis}.

Finally, for Theorem \ref{thm:c>0}, we need the asymptotic distribution of the ratio $R$ of the difference between the largest observed lifetime and the
largest uncensored lifetime taken as a proportion of
the largest lifetime.  For this, calculate, for $0<x<1$,
  \begin{align}\label{Rd}
    P[R>x]=&
P \Biggl[
\frac{
Y_u\Bigl(\frac{1}{1+\kappa} \Bigr) \vee Y_c\Bigl(\frac{\kappa}{1+\kappa}\Bigr)- Y_u\Bigl( \frac{1}{1+\kappa}\Bigr) }
{Y_u\Bigl(\frac{1}{1+\kappa} \Bigr) \vee Y_c\Bigl(\frac{\kappa}{1+\kappa}\Bigr)}    >x  \Biggr]  \cr
=&
\int_\R    P \Bigl[Y_c\Bigl( \frac{\kappa}{1+\kappa} \Bigr) \vee y
>\frac{y}{1-x}  \Bigr]  P \Bigl[Y_u\Bigl( \frac{1}{1+\kappa} \Bigr)\in \rmd y\Bigr]  \cr
=&
\int_\R (1-\Lambda^{\kappa/(1+\kappa)}(y/(1-x))\, 
\Lambda ^{1/(1+\kappa)} (\rmd y)\quad {\rm (using\ \eqref{yb})}  \cr
=&
\int_0^\infty \big(1- \exp\{-\kappa e^{-y/(1-x)}/(1+\kappa)\}
\big) \,  \rmd \big( \exp\{-e^{-y}/(1+\kappa)\}  \big).\cr
&
  \end{align}
  Some computations reduce the RHS of \eqref{Rd} to
  \be\label{fin}
      P[R>x]=
  \frac{1-x}{1+\kappa} 
  \int_0^\infty \big(1- e^{-\kappa u/(1+\kappa)}\big) \, 
  e^{-u^{1-x}/(1+\kappa)} u^{-x} \rmd u,
 \ee
 for $0<x<1$. When $x=0$ we can evaluate the integral and take the complement  to get  a mass at 0 for $R$ of 
 $1/(1+\kappa)$. 
Thus the cdf of $R$ can be written as the complement of \eqref{fin}.
With these we complete the proof of Theorem \ref{thm:c>0}. \halmos
%

%

\subsection{The case $\kappa=0$.}\label{subsec:c=0}
The case $\kappa=0$  of  very light censoring is of interest.
 When $\kappa=0$,
\eqref{e:relateFG} holds with $U$ slowly varying and from
\eqref{e:c0},
with $a(\cdot), b(\cdot)$ chosen as in \eqref{e:defab}, 
\begin{align*}
  P[  \max_{1\le i\le n}  U_{\Kgreater_i}    \leq a(n)x&+b(n)]
  =\Bigl( 1-\frac{nP[ U_{\Kgreater_1} > a(n)x+b(n)]}{n}\Bigr)^n\\
  =&\Bigl( 1-\frac{P[ U_{\Kgreater_1} > a(n)x+b(n)]}{P[ T_{\Kless_1} >
     a(n)x+b(n)]    }
\frac{nP[ T_{\Kless_1} >
     a(n)x+b(n)]    }{n}  \Bigr)^n \\
  \to & e^{-0} =1, \ x\in\R.
\end{align*}
Since this is true for any $x\in \R$, we have
$$\frac{ \max_{1\le i\le n}  U_{\Kgreater_i}  -b(n)}{a(n)} \Rightarrow -\infty.$$
So the analogue of \eqref{e:megilla} has limit
$(Y_u(1), -\infty)$ and \eqref{e:L1L2} has the degenerate limit
$$(L_1,L_1)=(Y_u(1),Y_u(1)). $$ Thus, for the case $\kappa=0$,
$$L=L_2-L_1=0$$
and as $n\to\infty$, the difference between the maximal observation
and the maximal uncensored lifetime vanishes asymptotically.

Consequently, in this situation we observe a vanishingly small level stretch at the righthand end of the KME, asymptotically. 
Likewise, $R=0$ w.p.1 in this case.


\section[Number]{Proof of Theorem \ref{th:limCts}}\label{sec:count}

We employ the {\it decoupage}  again, this time applying it
twice. Now the pairs $\{(T_i^*,U_i), i\geq 1\}$ are split  into independent
sets according to whether they are above the diagonal in the
$(t,u)$-plane or below.  Uncensored lifetimes $\{T^*_{\Kless _j }, 1\leq
j \leq N_u(n)\}$ depend on observations above the diagonal and the
independent random variable $N_u(n)$ and are independent of censored
observations which are  below the diagonal. The maximal
uncensored lifetime given in \eqref{mu},  $M_u(n)=\vee_{i=1}^{N_u(n)
}T_{\Kless_i}$, is independent of censored observations.

Take the i.i.d. collection $\{C_j:=U_{\Kgreater_j}, j\geq 1\}$ of
censored lifetimes. Some of these are (typically) less than $M_u(n)$ and some are greater.
The variable $M_u(n)$
is independent of censored observations.
We  condition on $M_u(n)=t$ and
split $\{C_j\}$ into two independent subsets
via the {\it decoupage\/},  into  the i.i.d. collection of those observations less than $t$ and those greater than $t$. 
The number of censored lifetimes
in the sample of size $n$ of lifetimes is  $N_c(n)=n-N_u(n)\sim np_c$.
Recall the definition of $N_c(>t)$ in \eqref{e:cs}.
(Jointly) conditional on $M_u(n)=t$ and $N_c(n)$, 
 $N_c(>t)$ 
is binomial with the number of trials being $N_c(n)$ and success
probability according to \eqref{e:tail>} being
\be\label{e:successprob}
p(t):=\frac{\int_t^\infty \bar F(s) G(\rmd s)}{\int_0^\infty \bar F(s)     G(\rmd s)}.
\ee
Throughout this section  we assume  \eqref{e:barF}, \eqref{e:barG} and \eqref{e:assume}   and keep $0<\kappa <\infty$.

\subsection{Asymptotics of $p(t)$.}\label{subsec:p(t)}
Since the argument requires that we condition on
$M_u(n)=t$, we consider the large sample distribution of $p(M_u(n))$.

\begin{proposition}\label{p2}
 We have with $p_c$ given
  in \eqref{e:pc},
  \beq \label{e:np(t)}
  np\bigl(M_u(n)\bigr) \Rightarrow \frac{\kappa}{p_c} E,
  \eeq
where $E$ is a standard exponential random variable.
\end{proposition}

\noindent{\bf Proof of Proposition \ref{p2}:}\ 
  Set
$Y_n:=(M_u(n)-b(n))/a(n)$ so $M_u(n)=a(n)Y_n+b(n).$ Recall from
\eqref{mu} and \eqref{e:L1L2} that
$$Y_n \Rightarrow Y_u(\frac{1}{1+\kappa}) \stackrel{d}{=} Y_u(1) +\log
\frac{1}{1+\kappa}.$$
By the Skorohod embedding theorem 
 (\cite[p.6]{resnickbook:2008})
convergence in distribution may be replaced by almost sure convergence. Doing this, we get
\begin{align*}
np(M_u(n))=& n   \frac{\int_{M_u(n)}^\infty \bar F(s)
             G(\rmd s)}{\int_0^\infty \bar F(s)     G(\rmd s)}&&{}\\
  \sim & \frac{n\kappa}{(1+\kappa)p_c} \bar H(M_u(n))
  \qquad (\text{from  \eqref{e:cenTail}})\\
= & n \bar H(b(n)) \frac{\kappa}{(1+\kappa)p_c}  \frac{\bar
    H(M_u(n))}{\bar H(b(n))}&& {}\\
  \sim & \frac{\kappa}{(1+\kappa)p_c}  \frac{\bar
         H(a(n)Y_n+b(n))}{\bar H(b(n))}&& {}\\
  \to & \frac{\kappa}{p_c} e^{-Y_u(1)}\stackrel{d}{=}
        \frac{\kappa}{p_c}E \qquad (\text{from \eqref{e:L1L2}}),
        \end{align*}
completing the proof of Proposition \ref{p2}.
  \halmos

  \subsection{Asymptotics of $N_c(>t)$.}\label{subsec:limCts}
  We first prove the  conditioned limit result in
Theorem \ref{th:limCts} and then remove the
  conditioning.
  Let
$$
\EE^{(n)}(\cdot)=\EE\bigl((\cdot)\,|\, M_u(n),\, N_c(n)\bigr)$$
be the conditional
expectation given $M_u(n)$ and $ N_c(n)$.


For $0<s<1$ the conditional generating function of the binomial rv
$N_c(>M_u(n))$ is
  \begin{align*}
    \EE^{(n)}\big( s^{N_c(>M_u(n))}    \big)
    =& \Bigl(1-p(M_u(n))(1-s)\Bigr)^{N_c(n)}\\
    =&  \Bigl(1-\frac{np(M_u(n))(1-s)}{n}\Bigr)^{n(N_c(n)/n)}.
    \end{align*}
Applying \eqref{e:np(t)} this converges to
\begin{align*}
 & \exp\{-\frac{\kappa}{p_c}E(1-s)p_c \}
    = \exp\{-\kappa E(1-s) \},
  \end{align*}
  which is the generating function of a Poisson random variable with
  parameter $\kappa E$.

  For the unconditional generating function, by dominated convergence
  \begin{align*}
    \EE s^{N_c(>M_u(n))}
    =&
\EE \Bigl(\EE^{(n)} s^{N_c(>M_u(n))}\Bigr) \\
    \to & 
  \EE\Bigl( \exp\{-\kappa E(1-s)\}\Bigr) \\
  =&
  \frac{1}{1+\kappa (1-s)} \\
  =&
  \frac{1-\frac{\kappa}{1+\kappa}}{1-s\frac{\kappa}{1+\kappa}},
  \end{align*}
  which is the generating function of the   geometric distribution.
  \halmos

\bibliography{bibfile.bib}

\def\cprime{$'$}
\begin{thebibliography}{10}

\bibitem{av:2018}
M.~Amica and I.~Van~Keilegom.
\newblock Cure models in survival analysis.
\newblock {\em Annual Review of Statistics and Its Application}, 5:311--342,
  2018.

\bibitem{bingham:goldie:teugels:1987}
N.H. Bingham, C.M. Goldie, and J.L. Teugels.
\newblock {\em Regular Variation}.
\newblock Cambridge University Press, 1987.

\bibitem{dehaan:1974equiv}
L.~de~Haan.
\newblock {Equivalence classes of regularly varying functions.}
\newblock {\em Stochastic Processes Appl.}, 2:243--259, 1974.

\bibitem{dehaan:ferreira:2006}
L.~de~Haan and A.~Ferreira.
\newblock {\em Extreme Value Theory: An Introduction}.
\newblock Springer-Verlag, New York, 2006.

\bibitem{emvz:2020}
M.~Escobar-Bach, R.A. Maller, I.~Van~Keilegom, and M.~Zhao.
\newblock Estimation of the cure rate for distributions in the {G}umbel maximum
  domain of attraction under insufficient follow-up.
\newblock {\em Preprint}, 2020.

\bibitem{ev:2018}
M.~Escobar-Bach and I.~Van~Keilegom.
\newblock Non-parametric cure rate estimation under insufficient follow-up
  using extremes.
\newblock {\em Journal of the Royal Statistical Society. Ser. B
  (Methodological)}, 2018.

\bibitem{gehan:1965}
E.~Gehan.
\newblock A generalized {W}ilcoxon test for comparing arbitrarily single
  censored samples.
\newblock {\em Biometrika}, 52:203--223, 1965.

\bibitem{maller:zhou:1993}
R.A. Maller and S.~Zhou.
\newblock The probability that the largest observation is censored.
\newblock {\em Journal of Applied Probability}, 30:602--615, 1993.

\bibitem{maller:zhoubook:1996}
R.A. Maller and Z.~Zhou.
\newblock {\em Survival Analysis with Longterm Survivors}.
\newblock Wiley, Chichester, first edition, 1996.

\bibitem{oblc:2012}
M.~Othus, B.~Barlogie, M.L. LeBlanc, and J.J. Crowley.
\newblock Cure models as a useful statistical tool for analyzing survival.
\newblock {\em Clinical Cancer Research}, 18:311--342, 2012.

\bibitem{pt:2014}
Y.~Peng and J.M.G. Taylor.
\newblock Cure models.
\newblock {\em In: Klein, J., van Houwelingen, H., Ibrahim, J. G., and Scheike,
  T. H., editors, Handbook of Survival Analysis, Handbooks of Modern
  Statistical Methods series, chapter 6. Chapman \& Hall, Boca Raton, FL,
  USA.}, pages 113--134, 2014.

\bibitem{resnickbook:2008}
S.I. Resnick.
\newblock {\em Extreme Values, Regular Variation and Point Processes}.
\newblock Springer, New York, 2008.
\newblock Reprint of the 1987 original.

\bibitem{ti:2014}
F.~Taweab and N.~A. Ibrahim.
\newblock Cure rate models: a review of recent progress with a study of
  change-point cure models when cured is partially known.
\newblock {\em J. Appl. Sci.}, 14:609--616, 2014.

\end{thebibliography}


\end{document}